\documentclass[12pt]{article}
\newtheorem{theor}{Theorem}
\newtheorem{prop}{Proposition}
\newtheorem{lem}{Lemma}

\begin{document}

\title{Once more about Voronoi's conjecture and space tiling zonotopes}

\author{
Michel Deza\footnote{Michel.Deza@ens.fr}
\\
Ecole Normale Sup\'erieure, Paris, and ISM, Tokyo
\and
Viacheslav Grishukhin\footnote{grishuhn@cemi.rssi.ru}
\\
CEMI, Russian Academy of Sciences, Moscow}

\date{}

\maketitle

\begin{abstract}
Voronoi conjectured that any parallelotope is affinely equivalent to a 
Voronoi polytope. A parallelotope is defined by a set of $m$ facet vectors 
$p_i$ and defines a set of $m$ lattice vectors $t_i$, $1\le i\le m$. We show 
that Voronoi's conjecture is true for an $n$-dimensional parallelotope $P$ if 
and only if there exist scalars $\gamma_i$ and  a positive definite 
$n\times n$ matrix $Q$ such that $\gamma_i p_i=Qt_i$ for all $i$. In this 
case the quadratic form $f(x)=x^TQx$ is the metric form of $P$. 

As an example, we consider in detail the case of a zonotopal parallelotope.
We show that $Q=(Z_{\beta}Z^T_{\beta})^{-1}$ for a zonotopal parallelotope 
$P(Z)$ which is the Minkowski sum of column vectors $z_j$ of the $n\times r$ 
matrix $Z$. Columns of the matrix $Z_{\beta}$ are the vectors 
$\sqrt{2\beta_j}z_j$, where the scalars $\beta_j$, $1\le j\le r$, are such 
that the system of vectors $\{\beta_jz_j:1\le j\le r\}$ is unimodular. $P(Z)$ 
defines a dicing lattice which is the set of intersection points of the 
dicing family of hyperplanes $H(j,k)=\{x:x^T(\beta_jQz_j)=k\}$, where $k$ 
takes all integer values and $1\le j\le r$.       

\end{abstract}



\section{Introduction}
A {\em parallelotope} is a convex polytope which fills the space facet to
facet by its translation copies without intersecting by inner points. 
Such a filling by parallelotopes is a {\em tiling}. Any parallelotope 
$P$ of dimension $n$ has the following three properties: 

{\em (symP) $P$ is centrally symmetric; 

(symF) each facet of $P$ is centrally symmetric; 

(proj) the projection of $P$ along any $(n-2)$-face is either a 
parallelogram or a centrally symmetric hexagon.} 

Venkov \cite{Ve} (and, independently, McMullen \cite{McM1}) proved that 
the above three properties are sufficient for a polytope to be a 
parallelotope. Aleksandrov \cite{Al}, knowing the Venkov's result and using 
his main idea, simplified the proof of Venkov and generalized his result. 

Coxeter \cite{Cox} noted that the condition $(proj)$ is necessarily and 
sufficient for a zonotope to be a parallelotope. But he considered only 
dimensions $n\le 3$. Shephard \cite{Sh} proved this assertion for $n=4$. 
Coxeter remarked also in \cite{Cox} that the condition $(proj)$ implies
that, for $n\le 3$, {\em every} parallelotope is a zonotope.   

A parallelotope is called $k$-{\em primitive} if each its $k$-face
(i.e. $k$-dimensional face) belongs to a minimal possible number 
$n-k+1$ of parallelotopes of its tiling. Besides, a $k$-primitivity 
implies the ($k+1$)-primitivity. A 0-primitive parallelotope is 
simply called {\em primitive}. Obviously, any parallelotope is 
$(n-1)$-primitive. 

The definition of $(n-2)$-primitivity shows that each $(n-2)$-face 
of a parallelotope belongs to at least 3 parallelotopes of its tiling. 
The property $(proj)$ above says that each $(n-2)$-face belongs to at 
most 4 parallelotopes. Hence the property $(proj)$ can be reformulated 
in the following form:

{\em (fhyp) each $(n-2)$-face of a parallelotope $P$ is contained in 
at most 3 affine hyperplanes supporting $(n-1)$-faces of its tiling.}

A {\em zonotope} $P(Z)$ is the Minkowski sum of a set $Z\subset 
{\bf R}^n$ of $n$-vectors. Every zonotope and all its faces are centrally 
symmetric. Hence, every zonotope satisfies the properties $(symP)$ and 
$(symF)$ of a parallelotopes. This implies that, for every $k$, 
$1\le k\le n-1$, the family of $k$-faces of $P(Z)$ is partitioned into 
{\em zones} of mutually parallel $k$-faces. For a zonotope $P(Z)$, every 
its zone of $k$-faces is in a one-to-one correspondence with a $k$-flat 
of the matroid $M(Z)$ represented by the set $Z$. Hence for a zonotope 
$P(Z)$ the property $(fhyp)$ is equivalent to the following property 
of the matroid $M(Z)$:

{\em (bin) every $(n-2)$-flat of $M(Z)$ is contained in at most three 
$(n-1)$-flats.}

This condition on a matroid is equivalent to its binarity. Since 
$M(Z)$ is trivially represented over the real field $\bf R$, $M(Z)$ 
is regular. This implies that 

{\em (mzr) a zonotope $P(Z)$ is a parallelotope if and only if the 
matroid $M(Z)$ is regular.} 

The regularity of $M(Z)$ implies that, for $z\in Z$, there are positive 
scalar $\beta_z$ such that the system of vectors $\{\beta_zz:z\in Z\}$ is 
unimodular. 
%The similarity of conditions $(iii)'$ and $(bin)$ was noted by 
%Jaeger in \cite{Ja}. He used $(bin)$ for to prove $(mrz)$.   

There is a special well known case of a parallelotope, namely, the 
Voronoi polytope $P_V(t_0)$ related to a point $t_0$ of a lattice $L$. 
It is the set of all points of the space that are at least as near 
to $t_0$ as to any other lattice point. Here the distance between two 
points $t$ and $t_0$ is the Euclidean norm $(t-t_0)^2$ of the vector 
$t-t_0$. Using an arbitrary positive definite quadratic form $f$ as 
a new norm $f(t-t_0)$ of the vector $t-t_0$, we obtain a parallelotope 
$P_f$, which we call the {\em Voronoi polytope with respect to the form} 
$f$. But any quadratic form $f(x)$ is affinely equivalent to the Euclidean 
form $x^2$. The corresponding affine transformation maps $P_f$ into $P_V$.  

The famous conjecture of Voronoi is that 

{\em (Voc) any parallelotope $P$ is affinely equivalent to a Voronoi 
polytope $P_V$.} 

A parallelotope is defined by a set of $m$ facet vectors $p_i$ and defines 
a set of $m$ lattice vectors $t_i$, $1\le i\le m$ (see (\ref{par}) below). 
We show that Voronoi's conjecture for an $n$-dimensional parallelotope $P$ 
is equivalent to the following condition

{\em (pQt) there exist scalars $\gamma_i$ and  a positive definite 
$n\times n$ matrix $Q$ such that $\gamma_i p_i=Qt_i$ for all $i$.} 

In this case the quadratic form $f(x)=x^TQx$ is the metric form of $P$. 

Giving an explicit matrix $Q$, Voronoi proved his conjecture for primitive 
parallelotopes. Since a primitive parallelotope $P$ is $(n-2)$-primitive, 
the projection of $P$ along any $(n-2)$-face gives a centrally symmetric 
hexagon. On the other hand, if any such projection gives a hexagon, then 
the parallelotope is $(n-2)$-primitive. Zhitomirskii \cite{Zh} remarked 
that the Voronoi's proof in \cite{Vo} uses only $(n-2)$-primitivity. This 
allowed him to extend the result of Voronoi for $(n-2)$-primitive 
parallelotopes. 

At this time all parallelotopes of dimension $n\le 5$ are known, and, since 
each of these parallelotopes is affinely equivalent to a Voronoi polytope, 
the Voronoi's conjecture is true for $n\le 5$. 

The two parallelotopes in ${\bf R}^2$, namely a centrally symmetric hexagon 
(primitive) and a parallelogram (non-primitive), are known since 
antiquity.   

In 1885, Fedorov described all 5 parallelohedra of dimension 3 including one 
primitive. 

Delaunay \cite{De} enumerated 51 four-dimensional parallelotopes (including 
known to Voronoi 3 primitive ones), and Shtogrin \cite{St} found the last 
missed by Delaunay 52th non-primitive parallelotope. 

Ryshkov and Baranovskii \cite{RB} found 221 primitive 5-dimensional 
parallelotopes. The missed 222th one is identified in \cite{EG}. Engel in 
\cite{En} and \cite{En1}, using a computer, enumerated 179 372 
parallelotopes of dimension 5 (including 222 primitive). 

Engel's proof of Voronoi's conjecture is based on a program which 
verified whether a given polytope is a parallelotope or not, and if it is, 
this program find an affinely equivalent Voronoi polytope. Similarly to a 
zonotope, a parallelotope has zones of mutually parallel edges. Such edge 
zone is called {\em closed} if any two-dimensional face has two or none 
edges of the zone. Otherwise, the zone is called {\em open}. Each 
closed zone can be contracted to an open zone. Conversely, some open 
zones can be extended to closed zones. Hence parallelotopes are 
partially ordered by inclusion of sets of closed zones. Maximal members 
of this order are primitive parallelotopes. So, any parallelotope can be 
obtained from primitive one by contraction of some closed zones. Note 
that a contracted Voronoi polytope is a parallelotope but, in general, 
it  is not a Voronoi polytope. By such a manner, Engel enumerated all 
5-dimensional parallelotopes. Giving a Voronoi polytope affinely 
equivalent to each parallelotope, Engel proved Voronoi's conjecture for 
dimension 5.     

McMullen \cite{McM} proved Voronoi's conjecture for parallelotopes 
which are zonotopes. Erdahl \cite{Er} gave another proof of the result 
of McMullen using a notion of a {\em lattice dicing} (which forms a 
lattice whose Voronoi polytope is a zonotope). 

We show that a zonotopal parallelotope $P(Z)$ is a Voronoi polytope 
with respect to the quadratic form $f(x)=x^TQx$, for 
$Q=(Z_{\beta}Z^T_{\beta})^{-1}$, where the columns of the matrix 
$Z_{\beta}$ are the vectors $\sqrt{2\beta_z}$ for $z\in Z$. Hence the 
transformation $x\to Ax$ maps $P(Z)$ into a Voronoi polytope if the 
matrix $A$ is a solution of the equation 
$A^TA=(2Z_{\beta}Z_{\beta}^T)^{-1}$.

The lattice of $P(Z)$ is a {\em dicing lattice} obtained as the set 
of intersection points of the family of hyperplanes $\{x:d^Tx=m\}$, 
where $m$ is an integer and the family of vectors $d=\beta_zQz$, 
$z\in Z$, being a linear transform of the unimodular system 
$\{\beta_zz:z\in Z\}$, forms also a unimodular system. 

\section{Parallelotopes}
Being a centrally symmetric polytope, any parallelotope $P$ with $m$ 
pairs of opposite facets has the following description.   
\begin{equation}
\label{par}
P=P(0)\equiv\{x\in {\bf R}^n: 
-\frac{1}{2}p_j^Tt_j\le p_j^Tx\le \frac{1}{2}p_j^Tt_j,\mbox{  }
1\le j\le m\}, 
\end{equation}
where $\frac{1}{2}t_j$ is the center of the facet 
\[F_j=\{x\in P(0):p_j^Tx=\frac{1}{2}p_j^Tt_j\} \]
and $p_j$ is a {\em facet vector} of the facet $F_j$. Obviously, $p_j$ 
is determined up to a scalar multiple. 

The parallelotope $P(t_j)$ is obtained by translation of $P(0)$ on the 
{\em translation vector} $t_j$. It is adjacent to $P(0)$ by the facet 
$F_j$. The set ${\cal T}=\{t_j:1\le j\le m\}$ of translation vectors 
generates a lattice $L({\cal T})$ whose points are the centers of 
parallelotopes of the tiling. We denote by $\pm\cal T$ the set of all 
translation vectors and their opposite. 

If, for all $j$, $\gamma_jp_j=t_j$ with scalar $\gamma_j>0$, then the 
parallelotope $P(0)$ is the Voronoi polytope $P_V(0)$ related to the 
zero point of the lattice $L({\cal T})$. 

Note that the usual Euclidean norm $x^2=x^Tx$ is used in the definition 
of the Voronoi polytope $P_V(0)$. But we can use an arbitrary positive 
quadratic form $f_Q(x)=x^TQx$ as a norm of $x$, where $Q$ is a symmetric 
positive definite $n\times n$ matrix. Hence we call the parallelotope 
\[P_Q(0)\equiv\{x\in {\bf R}^n: 
x^TQx\le (x-t)^TQ(x-t), \mbox{  }t\in\pm{\cal T} \}. \]
the {\em Voronoi polytope with respect to the form} $f_Q(x)$. 
Such a parallelotope relates to the lattice $L({\cal T})$. 

Now, if $Q=I_n$, i.e. $Q$ is the unit matrix, then $P_{I_n}(0)=P_V(0)$
is the Voronoi polytope of the lattice $L({\cal T})$. 

\begin{prop}
\label{PfV}
Let $P$ be a parallelotope, defined in (\ref{par}), $A$ be a 
non-singular $n\times n$ matrix such that $Q=A^TA$ is a 
symmetric positive definite matrix. Then the following assertions 
are equivalent:

(i) the map $x\to Ax$ transforms $P$ into a Voronoi polytope, i.e.,
Voronoi conjecture is true for $P$; 

(ii) $P$ is the Voronoi polytope $P_Q$ with respect to the positive 
definite quadratic form $f(x)=x^TQx$;

(iii) one can choose scalar $\gamma_j$ such that the following equalities 
hold  
\begin{equation}
\label{pQt}
\gamma_jp_j=Qt_j \mbox{ for all }1\le j\le m. 
\end{equation}
\end{prop}
{\bf Proof}. 
(i)$\Leftrightarrow$(ii)
Let the map $x\to x'=Ax$ maps $P$ into the Voronoi polytope 
\[P_V=\{x':{x'}^2\le (x'-t)^2, \mbox{ for }t\in\pm A{\cal T}\}.\] 
The converse transformation $x'\to x=A^{-1}x'$ maps $x'^2$ into 
\[(Ax)^TAx=x^TA^TAx=x^TQx.\]  
Hence the transformation $x'\to x$ maps $P_V$ into the Voronoi 
polytope $P_Q$ with respect to the form $f(x)=x^TQx$, which should 
coincide with $P$.  

(iii)$\Leftrightarrow$(ii)
Let $t_j$ and $\gamma_jp_j$ are connected by the relation (\ref{pQt}).
Then we can rewrite the inequality 
$\gamma_jp^T_jx\le \frac{1}{2}\gamma_jp^T_jt_j$, 
determining the facet $F_j$ of $P$, as $0\le t^T_jQt_j-2t_j^TQx$, or, 
adding $x^TQx$ to both the sides of this inequality, as 
\[x^TQx\le (x-t_j)^TQ(x-t_j). \]
In other words, we obtain that $P=P_Q$. This reasoning can be reversed. 
Hence we are done. 

%When $Q$ is symmetric and positive definite, then it can be decomposed  
%into a product $Q=A^TA$, where $A$ is a non-degenerate $n\times n$ 
%matrix. The quadratic form takes the form 
%\[f_Q(x)=(Ax)^T(Ax). \]
%We can consider $Ax$ as the image of $x$ under the affine transformation 
%$x\to Ax$. Then the Voronoi polytope $P_V(0)$ can be considered as an affine 
%image of the parallelotope $P_Q(0)$. $P_V(0)$ relates to the lattice  
%which is an affine image of the lattice of $P_Q(0)$. 
%The converse transformation $x\to A^{-1}x$ maps the polytope $P_Q(0)$ into 
%the Voronoi polytope $P_V(0)$. 

\section{Matroids}
Any set $X\subset {\bf R}^n$ of vectors gives a representation (or a 
coordinatization) (over the field ${\bf R}$) of a matroid $M=M(X)$. 
A matroid on a set $X$ is uniquely defined by its {\em rank function} 
rk$Y \le |Y|$ for all $Y\subseteq X$. Rank of a subset $Y\subseteq X$, 
is dimension of the space spanned by $Y$. A non-zero vector $x\in X$ 
has rk$x=1$ and is called a {\em point} of $M(X)$. Rank of the matroid 
$M(X)$ is the rank of $X$. A maximal by inclusion subset of rank $k$ is 
called a $k$-{\em flat}, or simply {\em flat}. So, a 1-flat is a point. 
$2$-flats are called {\em lines}. Let rk$M(X)=n$. Then $n-2$- and 
$n-1$-flats are called {\em colines} and {\em copoints}, respectively. 
So, a copoint $H$ spans a hyperplane $h$ in ${\bf R}^n$. For a detailed 
information on matroids, see, for example, \cite{Aig} and \cite{Wh}.    

For our aims, binary and regular matroids are important. A matroid is 
called {\em binary} if it is represented over a field of characteristic 
two. Another equivalent characterization of a binary matroid is $(bin)$ 
in terms of copoints and colines (see Theorem 7.22(iv) of \cite{Aig}):

\vspace{2mm} 
{\em (bim) A matroid is binary if and only if each its coline is contained 
in at most 3 copoints.} 

\vspace{2mm}
A binary matroid is called {\em regular} (or {\em unimodular}) if it is 
represented also over a field of characteristic distinct from two (see 
Theorem 7.35(iii) of \cite{Aig}). Any regular matroid $M(X)$ can be 
represented by a unimodular system of vectors $X$. A set $X$ of 
$n$-dimensional vectors is called {\em unimodular} if for {\em every} 
subset $B\subseteq X$ of $n$ independent vectors any vector of $X$ is 
represented as an integer linear combination of vectors from $B$. If we 
suppose that $X$ and $B$ are columns of the corresponding matrices, then 
the matrix $B^{-1}X$ is {\em totally unimodular}, i.e. all its minors are 
equals to 0 or $\pm 1$.

One of many equivalent characterizations of a unimodular system 
is as follows (cf. Theorem 7.37(v) of \cite{Aig}, Theorem 3.11(7) of 
\cite{Wh} and the condition V of \cite{McM}) (we denote by $y^Tx$ the 
scalar products of column vectors $x$ and $y$).

Let $H\subset X$ be a copoint of $M(X)$. Then $H$ generates a hyperplane 
$h\subset{\bf R}^n$. Let $p_H$ be a vector orthogonal to $h$. Obviously, 
the length of $p_H$ can be arbitrary.    

\vspace{2mm}
{\em (uni) A system of vectors $X$ is unimodular if and only if, for any 
copoint $H\subset X$ of the matroid $M(X)$, one can choose a scalar 
multiple $\gamma_H$ of the vector $p_H$ such that 
$\gamma_H p_H^Tx\in\{0,\pm 1\}$ for all $x\in X$.}

\vspace{2mm}
A proof of $(uni)$ is easily obtained from the fact that 
any vector of $X$ has $(0,\pm 1)$-coordinates in any base of $X$. 

Unfortunately, not every representation over $\bf R$ of a regular 
matroid is unimodular. But {\em every} representation $X$ over $\bf R$ 
can be transformed into a  unimodular representation by multiplying each 
vector $x\in X$ on an appropriate positive scalar multiple $\beta_x$. 
Hence $(uni)$ can be reformulated as follows. 

\vspace{2mm}
{\em (reg) A matroid $M(X)$ is regular if and only if one can choose scalar 
multiples $\gamma_H$ and $\beta_x$ of the vectors $p_H$ and $x$, 
respectively, such that $\gamma_H p_H^T\beta_x x\in\{0,\pm 1\}$ for any 
copoint $H\subset X$ of $M(X)$ and all $x\in X$.}

\section{Zonotopes}
A zonotope is the Minkowski sum of segments. If the zonotope is 
$n$-dimensional, then each segment $S_i$ is defined by an $n$-dimensional 
vector $z_i$ such that 
\[S_i=\{x\in {\bf R}^n:-z_i\le x\le z_i\}. \]
Let a zonotope $P(Z)$ be the Minkowskii sum of $r$ segments $S_i$ 
and be defined by a system of column vectors $\{z_i:1\le i\le r\}$ of an 
$n\times r$ matrix $Z$. Then 
\[P(Z)=\{x\in {\bf R}^n: x=Zy, \mbox{  }-1\le y_i\le 1, \mbox{  }
1\le i\le r\}. \]
Obviously, $P(Z)$ is centrally symmetric. Each family of all mutually 
parallel $k$-faces of $P(Z)$ is called a {\em zone}.

We can consider the matrix $Z$ as a set of column vectors and denote the 
set by the same letter $Z$. Let $M(Z)$ be the matroid represented by $Z$. 
Let $X\subseteq Z$ be a $k$-flat. Obviously, $X$ generates a sub-matrix of 
$Z$ of rank $k$. The set $X$ defines a zonotope $P(X)$. In particular, 
$P(\{z_i\})=S_i$. 

An important property connecting $P(Z)$ and $M(Z)$ is as follows 

\vspace{2mm}
{\em (face) Every $k$-face $F$ of $P(Z)$ has the form $F=F(X)$ with  
\[F(X)=\sum_{i:z_i\not\in X}\varepsilon_iz_i+P(X),\]
where $X$ is a $k$-flat of $M(Z)$ and $\varepsilon_i\in\{\pm 1\}$. 
Conversely, every $k$-flat $X$ of $M(Z)$ defines a $k$-face $F(X)$ of 
$P(Z)$ for some $\varepsilon_i\in\{\pm 1\}$. If $F(x)$ is a facet and 
$p$ is a vector normal to the hyperplane supporting $P(Z)$ in $F(X)$, 
then}   
\begin{equation}
\label{eps}
\varepsilon_i=\left\{\begin{array}{rl}
                   1, & \mbox{ if }p^Tz_i>0, \\ 
                  -1, & \mbox{ if }p^Tz_i<0, \\
                   0, & \mbox{ if }p^Tz_i=0. \\
                      \end{array}\right. 
\end{equation} 

\vspace{2mm}
This classical form of a face can be found in almost all papers on 
zonotopes. The fact that $X$ is a flat is explicitly given, for example, 
in Proposition 2.2.2 of \cite{B-Z}. In \cite{Ja} this fact is given in 
terms of a representation of $M(Z)$ by a chain group.   

We see that any face of a zonotope $P(Z)$ and $P(Z)$ itself are 
centrally symmetric. The center of a facet $F$ is given by 
\begin{equation}
\label{cnt}
\frac{1}{2}t=\sum_{i=1}^r\varepsilon_iz_i. 
\end{equation} 

Theorem~\ref{main1} below is well known. This is Proposition 3.3.4 of 
\cite{Wh} and Theorem 2.2.10 of \cite{B-Z}, both given there without 
proofs. We give a very short proof.  
\begin{theor}
\label{main1}
Let $Z$ be a set of vectors. The following assertions are equivalent:

(i) the zonotope $P(Z)$ is a parallelotope;

(ii) the matroid $M(Z)$ is regular. 
\end{theor}
{\bf Proof}. (i)$\Rightarrow$(ii) If $P(Z)$ is a parallelotope of 
dimension $n$, then its $(n-2)$-faces have the property $(fhyp)$. Since 
by $(face)$ each $k$-face of $P(Z)$ uniquely determines a $k$-flat of 
$M(Z)$, the matroid $M(Z)$ satisfies the condition $(bin)$. Hence $M(Z)$ 
is binary. Since $M(Z)$ is obviously represented over $\bf R$, by 
definition of a regular matroid, $M(Z)$ is regular. 

(ii)$\Rightarrow$(i) Since any zonotope is centrally symmetric and has 
centrally symmetric faces, $P(Z)$ has the first two properties $(symP)$ 
and $(symF)$ of a parallelotope. Hence we have to prove that $P(Z)$ has 
the property $(fhyp)$. Since $M(Z)$ is regular, it is binary and satisfies  
the condition $(bin)$. This condition and $(face)$ imply that $P(Z)$ has 
the property $(fhyp)$, too. 
  
\vspace{3mm}
Note that in \cite{McM} and \cite{Sh}, the condition $(fhyp)$ of a 
parallelotope is formulated as the condition $(bin)$ (see II of \cite{McM}
and (11) of \cite{Sh}) but without mention of a notion of matroid. 

The above proof of the implication (i)$\Rightarrow$(ii) can be found in 
Section 4.2 of \cite{Ja}. 
 
Now we show that a zonotopal parallelotope $P(Z)$ is a Voronoi polytope 
with respect to a quadratic form $f_D(x)$ for some matrix $D$. 

Let $\beta_i$ be the positive multiple of $z_i$ mentioned in the condition  
$(reg)$. So, the system $\{\beta_iz_i:1\le i\le r\}$ is 
unimodular. Let $Z_{\beta}$ be an $n\times r$ matrix whose columns are 
the vectors $\sqrt{2\beta_i}z_i$. 
\begin{lem}
\label{pt}
Let $P(Z)$ be a full-dimensional zonotopal parallelotope determined by 
a system $Z$. Let the facet vectors $p_j$, $1\le j\le m$, satisfy the 
condition $(reg)$ (with $\gamma_H=\gamma_j$). Then the 
following equalities hold
\[\gamma_jp_j=Qt_j, \mbox{  }1\le j\le m, 
\mbox{ where } Q=(Z_{\beta}Z_{\beta}^T)^{-1}. \]
\end{lem}
{\bf Proof}. By (\ref{cnt}), the center of a facet $F_j$ is   
\[\frac{1}{2}t_j=\sum_{i=1}^r\varepsilon_{ij}z_i, \]
where $\varepsilon_{ij}$ is $\varepsilon_i$ of (\ref{eps}) for $p=p_j$.      
If the facet vectors $p_j$, $1\le j\le m$, satisfy the condition  
$(reg)$, then $\varepsilon_{ij}=\gamma_j\beta_i p_j^Tz_i$. Hence we have
\[t_j=2\sum_{i=1}^r\beta_i\gamma_j(p_j^Tz_i)z_i=
\gamma_j\sum_{i=1}^r\sqrt{2\beta_i}z_i(\sqrt{2\beta_i}z_i^Tp_j)
=Z_{\beta}Z^T_{\beta}\gamma_jp_j.\]

Since $P(Z)$ is full-dimensional, the matrix $Z_{\beta}Z^T_{\beta}$ is 
a positive definite symmetric matrix. So, it is non-singular, and we 
have $\gamma_jp_j=Qt_j$. The result follows. 

\vspace{2mm}
Hence we have the following 
\begin{theor}
\label{main2}
The map $x\to Ax$ transforms a zonotopal parallelotope $P(Z)$ into a 
Voronoi polytope, where the matrix $A$ is a solution of the equation 
$A^TA=(Z_{\beta}Z^T_{\beta})^{-1}$.
\end{theor}    

\section{Dicings}

Let $D$ be a set of vectors spanning ${\bf R}^n$. Then the set $D$ defines 
a family ${\cal H}(D)$ of parallel hyperplanes 
$H(d,m)=\{x\in {\bf R}^n:x^Td=m\}$, where $d\in D$ and $m$ is an integer. 
Let $B\subseteq D$ be a basis of $D$. Then the set of intersection points 
of the hyperplanes $H(d,m)$ for $d\in B$ and all $m$ is a lattice $L(B)$. 
Erdahl and Ryshkov \cite{ER} proved that the set of intersection points of 
the hyperplanes of the whole family ${\cal H}(D)$ is a lattice $L(D)$ 
(which then coincides with $L(B)$) if and only if $D$ is a unimodular system. 
(In \cite{ER} the notion of a unimodular system is not mentioned). In this 
case, the family ${\cal H}(D)$ is called a {\em lattice dicing}. We call the 
lattice $L(D)$ a {\em dicing lattice}. The family ${\cal H}(D)$ dices 
${\bf R}^n$ into polytopes which are Delaunay polytopes of $L(D)$.  

Now we show that the lattice $L({\cal T})$ formed by the centers of 
zonotopal parallelotopes is a dicing lattice. 

Take the unimodular system $\{\beta_iz_i:1\le i\le r\}$ considered in the 
preceding sections. Obviously, the system $D=\{\beta_iQz_i:1\le i\le r\}$ 
is also unimodular. Hence it defines a dicing lattice $L(D)$. We show that 
$L(D)=L({\cal T})$, i.e. that $L(D)$ is the lattice of the zonotopal 
parallelotope $P(Z)$. 

Let $I_j=\{i:p_j^Tz_i=0\}$. So, $p_j$ is the intersection of the hyperplanes 
\[ H_i=\{x\in {\bf R}^n:x^Tz_i=0\},  \]
for $i\in I_j$. Similarly, $t_j$ lies in the intersection of the hyperplanes 
$G_i$, $i\in I_j$, of the family ${\cal H}(D)$, where 
\[ G_i=\{x\in {\bf R}^n:x^T\beta_iQz_i=0\}. \]
In fact, since $Qt_j=\gamma_jp_j$, we have, for $i\in I_j$,  
\[t_j^TQz_i=(Qt_j)^Tz_i=\gamma_jp^T_jz_i=0. \]
Since $Q$ is non-degenerate, dimension of the intersection of the hyperplanes  
$G_i$, $i\in I_j$, coincides with dimension of the interaction of the 
hyperplanes $H_i$, $i\in I_j$ which is equal 1. Hence $t_j$ spans the above 
intersection. Besides, if $t_k$ does not lie in $G_i$, we have 
$t^T_k\beta_iQz_i=\beta_i\gamma_kp_k^Tz_i=\varepsilon_{ik}\in\{\pm 1\}$. 
In other words, the endpoints of these $t_k$'s lie in the hyperplanes of 
the family $D$ parallel to $G_i$ and neighboring to it. 
This implies that $L({\cal T})=L(D)$.  

Similarly to the facet vectors $p_j$ determining facets of the parallelotope 
$P(Z)$, the dicing vectors $d_i$ determines facets of Delaunay polytopes. 
Edges of these Delaunay polytopes are just the lattice vectors $t_j$. 
 
Using the dicing vectors $d_i=\beta_iQz_i$, $1\le i\le r$, and the equality 
$Q^{-1}=Z_{\beta}Z^T_{\beta}=\sum_{i=1}^r2\beta_iz_iz^T_i$, 
we can represent the quadratic form $f(x)=x^TQx$ as a sum of quadratic forms 
of rank 1:
\[f(x)=x^TQx=x^TQQ^{-1}Qx=x^TQ(\sum_{i=1}^r2\beta_iz_iz^T_i)Qx=
\sum_{i=1}^r\frac{2}{\beta_i}(x^Td_id^T_ix)=
\sum_{i=1}^r\lambda_i(d_i^Tx)^2,\] 
where $\lambda_i=\frac{2}{\beta_i}$. 
 
Erdahl proves in \cite{Er} that a parallelotope related to a dicing lattice 
defined by dicing vectors $d\in D$ is a Voronoi polytope with respect 
to the quadratic form $f(x)=\sum_{d\in D}\lambda_d(d^Tx)^2$, $\lambda_d>0$, 
which is a zonotope $P(Z)$. He gives linear expressions connecting the 
vectors $z\in Z$ generating the zonotope $P(Z)$ and the vectors $d\in D$. 

%We saw that the quadratic form of a zonotope is uniquely represented as 
%a weighted sum of $r$ quadratic forms $(d_i^Tx)^2$ of rank 1. As the vectors 
%of a space of dimension $\frac{n(n+1)}{2}$, these functions are linearly 
%independent. Recall that quadratic forms of rank 1 are extreme rays of 
%L-type domains. Hence the lattice $L({\cal T})$ of a zonotopal parallelotope 
%belongs to a simplicial L-type domain of dimension $r$. Since {\em 
%non-rigidity degree} of a lattice is dimension of its L-type domain, we 
%have that non-rigidity degree of $L({\cal T})$ is equal to $r$, i.e. to 
%the number of zone vectors of the zonotope $P(Z)$.   

%\newpage

\end{document}